# Dodecahedral Structures with Mosseri-Sadoc Tiles

**Nazife Ozdes Koca[a,*], Ramazan Koc[b], Mehmet Koca[c], and Abeer Al-Siyabi[a]**

[a]Department of Physics, College of Science, Sultan Qaboos University, P.O. Box 36, Al-Khoud 123, Muscat, Sultanate of Oman, *Correspondence e-mail: nazife@squ.edu.om
[b]Department of Physics, Gaziantep University, Gaziantep, Turkey
[c]Department of Physics, Cukurova University, Adana, Turkey, retired

## ABSTRACT

3D-facets of the Delone cells of the root lattice $D_6$ which tile the six-dimensional Euclidean space in an alternating order are projected into three-dimensional space. They are classified into six Mosseri-Sadoc tetrahedral tiles of edge lengths 1 and golden ratio $\tau = \frac{1+\sqrt{5}}{2}$ with faces normal to the 5-fold and 3-fold axes. The icosahedron, dodecahedron and icosidodecahedron whose vertices are obtained from the fundamental weights of the icosahedral group are dissected in terms of six tetrahedra. A set of four tiles are composed out of six *fundamental* tiles, faces of which, are normal to the 5-fold axes of the icosahedral group. It is shown that the 3D-Euclidean space can be tiled face-to-face with maximal face coverage by the composite tiles with an inflation factor $\tau$ generated by an inflation matrix. We note that dodecahedra with edge lengths of 1 and $\tau$ naturally occur already in the second and third order of the inflations. The 3D patches displaying 5-fold, 3-fold and 2-fold symmetries are obtained in the inflated dodecahedral structures with edge lengths $\tau^n$ with $n \geq 3$. The planar tiling of the faces of the composite tiles follow the edge-to-edge matching of the Robinson triangles.

# 1. Introduction

After the discovery of the icosahedral quasicrystal by D. Shechtman (Shechtman *et al*, 1984) quasicrystallography gained enormous momentum. Its point symmetry can be described by the Coxeter group $H_3$ representing the icosahedral symmetry of order 120. Recent developments indicate that the quasicrystals exhibit 5-fold, 8-fold, 10-fold, 12-fold, and 18-fold symmetries. For a general exposition we refer the reader to the references on quasicrystallography (Di Vincenzo & Steinhardt, 1991; Janot, 1993; Senechal, 1995; Steurer, 2004; Tsai, 2008).

Aperiodic tilings in general (Baake & Grimm, 2013; Baake & Grimm, 2020) and the icosahedral quasicrystallography in particular constitute the main theme of research for many scientists from diverse fields of interest. The subject is mathematically intriguing as it requires the aperiodic tilings of the space by some prototiles. Projection technique from higher dimensional lattices is a promising approach. For an $(n + 1)$-fold symmetric planar aperiodic tilings one can use the projections of the 2D-facets of the Delone and Voronoi cells of the $A_n$ lattices (Koca et. al., 2019). Aperiodic tilings of three-dimensional Euclidean space with icosahedral symmetry have been a challenging problem. A tiling scheme with 7-prototiles has been proposed by Kramer (Kramer, 1982) which was converted to a four-tile model by Mosseri and Sadoc (Mosseri & Sadoc, 1982) leading to a 6-tetrahedral tiling system. Later it has been shown by Kramer and Papadopolos (Kramer & Papadopolos, 1994) that the tetrahedral tiles can be obtained from the root lattice $D_6$ by cut-and- project technique. See also the reference (Papadopolos & Ogievetski, 2000) for further information.

There have been two more approaches for the aperiodic order of the 3D Euclidean space with icosahedral symmetry. With an increasing order of symmetry, the set of four prototiles (Socolar & Steinhardt, 1986) consists of acute rhombohedron, Bilinski dodecahedron, rhombic icosahedron and rhombic triacontahedron. They are obtained from the Ammann tiles of acute and obtuse rhombohedra as they are the building blocks of the above composite tiles. A decoration scheme of the Ammann tiles were proposed by Katz (Katz, 1989) and it has been recently revived by Hann-Socolar-Steinhardt (Hann, Socolar & Steinhardt, 2018). Danzer (Danzer, 1989) proposed a more fundamental tiling scheme known as the *ABCK* tetrahedral prototiles with its octahedral tiling is denoted by $< ABCK >$. A common feature between the Amman rhombohedra (thereof their composite tiles) and the octahedral Danzer tiles is the fact that their faces are all normal to the 2-fold axes of the icosahedral group. Eventually, it was later shown that these two sets of tiles are related to each others (Danzer, Papadopolos & Talis, 1993; Roth, 1993). Danzer tiles are projected from $D_6$ lattice with the cut-and-project scheme (Kramer et. al., 1994a). Ammann rhombohedra can be obtained from the projections of the six-dimensional cubic lattice represented by the Coxeter-Dynkin diagram $B_6$ (Koca, Koca & Koc, 2015). Vertices of the Danzer prototiles can be derived from the fundamental weights of the icosahedral group which in turn can be obtained from the root lattice $D_6$ (Al-Siyabi, Koca & Koca, 2020). Kramer and Andrle (Kramer & Andrle, 2004) have also investigated the Danzer tiles in the context of $D_6$ lattice with its relation to the wavelets.

In this paper we demonstrate that the Mosseri-Sadoc tetrahedral tiling model can be obtained from the projections of the Delone cells of the root lattice $D_6$ without invoking the cut-and- project technique. Since the Delone cells tile the root lattice in an alternating order (Conway & Sloane, 1999) it is expected that the tiles projected from the Delone cells may tile the 3D-Euclidean space in an aperiodic manner with an icosahedral symmetry. See the reference (Koca et. al., 2018) for a detailed exposition of the root lattice $D_6$ discussing the Voronoi cell and the Delone polytopes. The dual of the root polytope of $D_6$ determined as the orbit of the weight vector $\omega_2$ is the Voronoi cell which is the disjoint union of the polytopes represented as the orbits of the weight vectors $\omega_1, \omega_5$ and $\omega_6$. These vectors represent the *holes* of the root lattice whose distance from the lattice is a local maximum. The distance of either $\omega_5$ or $\omega_6$ from the lattice is an *absolute maximum* and they are called the *deep holes* and the hole represented by the vector $\omega_1$ is called the *shallow hole* (Conway & Sloane, 1999, p. 33). The orbit of $\omega_1$ represents a cross polytope and either of the orbits of $\omega_5$ and $\omega_6$ is known as a hemicube or half measure polytope (Coxeter, 1973, p. 155) whose vertices are the alternating vertices of a 6D-cube. Delone cells are the polytopes whose vertices are the lattice vectors; Delone cells centralize the vertices of the Voronoi cell. For example, the Delone cells centralizing the vertices of the Voronoi cell $V(0)$ of the root lattice centered around the origin can be represented as the orbits of the weight vectors $\omega_1, \omega_5$ and $\omega_6$ shifted according the rule: $\omega_1 + \omega_1$, $\omega_5 + \omega_5$ and $\omega_6 + \omega_6$. Here, with $\omega_1 + \omega_1$, for example, we mean that one vector of the orbit of $\omega_1$ is added to all vectors of the orbit of $\omega_1$. With this we obtain 12 cross polytopes centralizing the vertices $\pm l_i (i = 1,2,...,6)$ of the Voronoi cell $V(0)$. Similar arguments are valid for the hemicubes.

When projected into 3D-Euclidean space the cross polytope turns into an icosahedron and each hemicube decomposes as the disjoint union of a dodecahedron and an icosahedron. The 240 tetrahedral facets of the Delone cell represented by the orbit of the weight vector $\omega_1$ projects into four types of tetrahedral tiles of edge lengths 1 and $\tau$. Similarly, 640+640 3D-facets of the hemicubes project into six tetrahedral tiles (including former four tiles) dissecting the dodecahedron and icosahedron. These are the Mosseri-Sadoc (Mosseri-Sadoc, 1982) tiles which we call them the *fundamental* tiles. The edge lengths of the equilateral triangular faces of the fundamental tiles being normal to the 3-fold axes are of two types 1 and $\tau$ and as such they cannot be partitioned in terms of each other. For this, we define a new set of four *composite* prototiles assembled by the fundamental tiles whose faces are normal to the 5-fold axes only and their faces consist of Robinson triangles. They are also defined in the reference (Mosseri-Sadoc, 1982) and studied in (Papadopolos & Ogievetsky, 2000) and (Kramer & Papadopolos, 1994b). The composite tiles can then be inflated by an inflation factor $\tau$ with an inflation matrix. In the following, the procedure is described as to how the 3D-space with four composite tiles is tiled and the emergence of the dodecahedral structures are discussed. This paper is an expanded version of the paper by Koca et. al. (Koca et. al., 2020) discussing the details of the dissections of the fundamental icosahedral polyhedra and displaying the 5-fold, 3-fold and 2-fold symmetries of the dodecahedral patches.

Before we proceed further, we should mention about the projection of the 4D and 5D-facets of the Delone cells either in the form of cross polytope or hemicube. Their facets are regular simplexes, such as, a 3D-facet is a tetrahedron, 4D-facet is a 4-simplex or sometimes called a 5-cell and so on. Projections of higher dimensional facets can be determined as the union of the lower dimensional facets. Since we will consider the projections of the polytopes in 6D-space into 3D-space all higher dimensional facets projects into polyhedra which are unions of the projected 3D-facets. Therefore, it is sufficient to work with the projections of the 3D-facets which are regular tetrahedra of the Delone cells of the $D_6$ lattice. Some examples of these will be given in Sec.2a.

The paper is organized as follows. In Sec. 2, we discuss the structures of the Delone cells of the root lattice $D_6$. Explicit dissections of the icosahedron and dodecahedron are studied in terms of the fundamental tiles. In Sec. 3 four composite Mosseri-Sadoc tiles are assembled with the fundamental tiles so that their faces consist of the Robinson triangles normal to the 5-fold axes. The composite tiles can be inflated with the $4 \times 4$ inflation matrix $M$ whose eigenvalues are $\tau^3$, $\tau, \sigma$ and $\sigma^3$ where $\sigma = -\tau^{-1} = \frac{1-\sqrt{5}}{2}$ is the algebraic conjugate of $\tau$. It is also pointed out that a simple modification of the composite tiles leads to similar dodecahedral structures. The right and left eigenvectors of the inflation matrix $M$ corresponding to the Perron-Frobenius (PF) eigenvalue $\tau^3$ are calculated and the projection matrix is formed as the tensor product of the right and left eigenvectors. The composite tiles display dodecahedral structures of edge lengths 1 and $\tau$ already in the 2$^{nd}$ and 3$^{rd}$ order of the inflation. In the increasing order of inflation, a systematic construction of the dodecahedral structures filling the space with the binding composite tiles are studied and the patches of dodecahedral structures are illustrated. The 5-fold, 3-fold and 2-fold planar symmetries of the patches are identified. In the conclusive remarks we propose two possible tiling schemes. Some details of the 3-fold symmetric construction of the dodecahedron $d(1)$ is described in Appendix A which also displays 5-fold symmetry. Appendix B describes another construction suitable for 5-fold symmetry only. Appendix C discusses the constructions of the icosidodecahedra of edge lengths 1 and $\tau$ denoted by $id(1)$ and $id(\tau)$ respectively.

## 2. Projections of the Delone cells of the root lattice $D_6$

We first discuss the projection of the Delone cell characterized as the orbit of the weight vector $\omega_1$ of the root lattice $D_6$, the point group of which, is of order $2^5 6!$. The orbit of $\omega_1$ known as the cross poytope in 6D space possesses 12 vertices, 60 edges, 160 triangular faces, 240 tetrahedral facets, 192 4-simplexes and 64 5-simplexes which can be directly obtained from the Coxeter-Dynkin diagram of $D_6$ in Fig. 1. More detailed discussion on the structure of the Delone cell related to the orbits of the weight vectors $\omega_5$ and $\omega_6$ will be given in Sec. 3b.

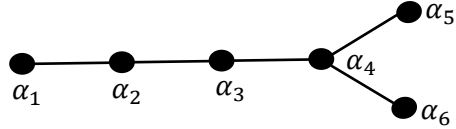

**Figure 1**
Coxeter-Dynkin diagram of the root lattice $D_6$.

For further use we will introduce the orthonormal set of vectors $l_i$, $(l_i, l_j) = \delta_{ij}$, $(i, j = 1, 2, \dots, 6)$ where the simple roots are given by $\alpha_i = l_i - l_{i+1}$, $(i = 1, 2, \dots, 5)$ and $\alpha_6 = l_5 + l_6$. The weight vectors $\omega_1$, $\omega_5$ and $\omega_6$ are given by the orthonormal set of vectors as

$$\omega_1 = l_1, \quad \omega_5 = \tfrac{1}{2}(l_1 + l_2 + l_3 + l_4 + l_5 - l_6), \quad \omega_6 = \tfrac{1}{2}(l_1 + l_2 + l_3 + l_4 + l_5 + l_6), \quad (2)$$

and their orbits under the point group $D_6$ are the sets of 12, 32 and 32 vectors respectively given by

$$\pm l_i, i = 1, 2, \dots, 6; \qquad (3a)$$

$$\tfrac{1}{2}(\pm l_1 \pm l_2 \pm l_3 \pm l_4 \pm l_5 \pm l_6). \qquad (3b)$$

Here odd and even combinations of the negative sign in (3b) correspond to the orbits of the weights $\omega_5$ and $\omega_6$ respectively. The set of 12 vectors in (3a) represents the vertices of a cross polytope and those 64 vectors in (3b) constitute the vertices of a 6D-cube and those vectors with either even or odd minus sign represent a hemicube in 6-dimensions. The vectors in (3a-b) also represent the vertices of the Voronoi cell $V(0)$ of the root lattice. A suitable representation of the orthonormal set of vectors $l_i$ in 6D-Euclidean space can be obtained as

$$\begin{bmatrix} l_1 \\ l_2 \\ l_3 \\ l_4 \\ l_5 \\ l_6 \end{bmatrix} =: \sqrt{\frac{2}{2+\tau}} \frac{1}{2} \begin{bmatrix} 1 & \tau & 0 & \tau & -1 & 0 \\ -1 & \tau & 0 & -\tau & -1 & 0 \\ 0 & 1 & \tau & 0 & \tau & -1 \\ 0 & 1 & -\tau & 0 & \tau & 1 \\ \tau & 0 & 1 & -1 & 0 & \tau \\ -\tau & 0 & 1 & 1 & 0 & \tau \end{bmatrix}, \qquad (4)$$

where the first three and last three components represent the vectors in the complementary $E_\parallel$ and $E_\perp$ spaces respectively. In what follows we shall represent the vectors $l_i$ by their first three components in the space $E_\parallel$ by deleting the overall factor $\sqrt{\frac{2}{2+\tau}}$ as

$$l_1 = \frac{1}{2}(1, \tau, 0), \quad l_2 = \frac{1}{2}(-1, \tau, 0), \quad l_3 = \frac{1}{2}(0, 1, \tau),$$
$$l_4 = \frac{1}{2}(0, 1, -\tau), \quad l_5 = \frac{1}{2}(\tau, 0, 1), \quad l_6 = \frac{1}{2}(-\tau, 0, 1). \qquad (5)$$

Note that we keep the same notation for the vectors $l_i$ to avoid the frequent use of the notation $\parallel$. This set of vectors are useful because they are directly related to the coordinates of the Danzer *ABCK* tiles (Al-Siyabi, Koca & Koca, 2020).

**a) Projection of the cross polytope into an icosahedron and its dissection**

The set of vectors $\pm l_i$, $(i = 1, 2, \dots, 6)$ in (4) representing a cross polytope in 6-dimensions now represents an icosahedron with 12 vertices given in (5), a sketch of which is given in Fig. 2.

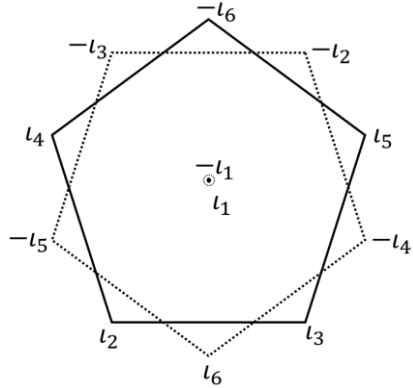

**Figure 2**
A sketch of icosahedron with vectors $\pm l_i$ of (5).

The 3D-facets of the cross polytope whose vertices are given by $\pm l_i$ of (4) are regular tetrahedra in 6D-space and are represented by four orthonormal set of vectors chosen from $\pm l_i$. For example, the set of unit vectors $l_1, l_2, l_3$ and $l_4$ given in (4) represents the vertices of a tetrahedron of edge length $\sqrt{2}$ which projects into a tetrahedron denoted by $t_1$ with 5 edge lengths 1 and one edge length $\tau$ as shown in Table 1 which displays six tetrahedral tiles projected from all Delone cells. Faces of six projected tetrahedral tiles consist of equilateral triangles of edges either 1 or $\tau$ and two types of Robinson triangles. To give another example, let us consider the regular tetrahedron in 6D-space with vertices $l_2, -l_4, -l_5$ and $-l_6$ which projects into a tetrahedron denoted by $t_6$ with 5 edge lengths $\tau$ and one edge length 1. It is clear from these two examples that the number of tetrahedra is the product $15 \times 2^4 = 240$, corresponding to the 4 choices out of 6 times the number of sign changes. The number of 3D-facets is also equal to the number of cosets of the group leaving the tetrahedron intact in the point group of order $2^5 6!$. One can check that 240 regular tetrahedra in 6D-space either project into one of those four tetrahedra $t_1, t_2, t_5$ and $t_6$ or isosceles trapezoids of edge lengths $(1,1,1,\tau)$ as seen from Table 1. They are all present in the icosahedron as 3-dimensional tiles or planar sections. Consequently, projection of the cross polytope leads to a dissection of the icosahedron with possible combinations of the tetrahedra $t_1, t_2, t_5$ and $t_6$ that will be explained below.

Perhaps it is more informative to give examples regarding the projections of the 4D and 5D-facets into 3D-space. The set of 5 vertices $l_1, l_2, l_3, l_4$ and $l_5$ in 6D-space constitute a 4-simplex consisting of 5 tetrahedral facets. When it is projected into 3D-space it leads to a polyhedron with 5 vertices, 8 edges and 5 faces. One of the edges is of length $\tau$ and the others are all of length 1. Three of the faces are equilateral triangles of edge lengths 1 and the other two faces are of types $(1,1,1,\tau)$ and $(1,1,\tau)$ which can be obtained as the composite tile $t_1 + t_2$. If one considers the projection of a 5-simplex with 6 vertices $l_1, l_2, l_3, l_4, l_5$ and $-l_6$ we obtain a pentagonal pyramid which is the union of the tetrahedral facets $t_1 + t_2 + t_1$.

An icosahedron can be constructed as the union of one pentagonal antiprism and two pentagonal pyramids or as the union of one Johnson solid $J_{63}$ and three pentagonal pyramids. Johnson solids are strictly convex polyhedra consisting of regular polygons (Johnson, 1966). We will consider the second choice for an illustration of the dissection of an icosahedron in terms of the four tetrahedra by following even a simpler method. We assemble the tiles $t_5$ and $t_6$ into a composite tile as

$$T_3 =: t_5 + t_6 + t_5 \tag{6}$$

by matching their equilateral triangular faces of edge length $\tau$ where $t_6$ is placed between two $t_5$ tiles. They form the "tent" of Kramer (Kramer, 1982) with $N_0 = 6$ vertices, $N_1 = 10$ edges and $N_2 = 6$ faces (henceforth $N_i$ will be used for vertices, edges and faces in this order) consisting of pentagonal base of edge length 1 and five Robinson triangles with edge lengths $(1, \tau, \tau)$. The 5 copies of tile $t_2$ with the faces $(1, \tau, \tau)$ can be matched face-to-face with the Robinson triangles of $T_3$, then the gaps between the successive tiles $t_2$ can be filled with 5 copies of $t_1$ on the faces $(1, 1, \tau)$. It is almost done except covering the pentagonal base by a pentagonal pyramid formed by the sandwich of tiles $t_1 + t_2 + t_1$. The result is an icosahedron $i(1)$ of edge length 1. It can be written as the union of the tiles $t_1, t_2, t_5$ and $t_6$ as

$$i(1) = 7t_1 + 6t_2 + 2t_5 + t_6 = 7t_1 + 6t_2 + T_3. \tag{7}$$

One possible set of coordinates of the tiles constituting the icosahedron in (7) can be represented by the set of vectors of (5) as

$t_5: (l_1, l_6, -l_3, -l_5);\quad t_6: (l_1, l_6, -l_2, -l_3);\quad t_5: (l_1, l_6, -l_2, -l_4);$
$t_2: (l_1, -l_6, -l_2, -l_3); (l_1, l_4, -l_3, -l_5); (l_1, l_2, l_6, -l_5); (l_1, l_3, l_6, -l_4); (l_1, l_5, -l_2, -l_4);$
$t_1: (l_1, l_4, -l_6, -l_3); (l_1, l_2, l_4, -l_5); (l_1, l_2, l_3, l_6); (l_1, l_3, l_5, -l_4); (l_1, l_5, -l_2, -l_6);$
$t_1: (-l_1, -l_3, -l_5, l_6); t_2: (-l_1, -l_2, -l_3, l_6); t_1: (-l_1, -l_2, -l_4, l_6).$ (8)

**Table 1**
The fundamental tiles projected from Delone cells of $D_6$ (All triangular faces with edge lengths $(1, 1, \tau)$ and $(1, \tau, \tau)$ are orthogonal to 5-fold axes; the faces with edge lengths $(1, 1, 1)$ and $(\tau, \tau, \tau)$ are all orthogonal to 3-fold axes).

| Name of tile | Sketch | Number of faces $(a, b, c)$ | Volume | Mosseri Sadoc notation | Kramer's notation |
|---|---|---|---|---|---|
| $t_1$ | | $2 \times (1, 1, 1)$ <br> $2 \times (1, 1, \tau)$ | $\dfrac{1}{12}$ | B | $B_\parallel^*$ |
| $t_2$ | | $1 \times (1, 1, 1)$ <br> $2 \times (1, 1, \tau)$ <br> $1 \times (1, \tau, \tau)$ | $\dfrac{\tau}{12}$ | G | $D_\parallel^*$ |
| $t_3$ | | $1 \times (\tau, \tau, \tau)$ <br> $3 \times (1, 1, \tau)$ | $\dfrac{\tau}{12}$ | E | $G_\parallel^*$ |
| $t_4$ | | $1 \times (1, 1, 1)$ <br> $3 \times (1, \tau, \tau)$ | $\dfrac{\tau^2}{12}$ | F | $F_\parallel^*$ |
| $t_5$ | | $1 \times (\tau, \tau, \tau)$ <br> $1 \times (1, 1, \tau)$ <br> $2 \times (1, \tau, \tau)$ | $\dfrac{\tau^2}{12}$ | C | $C_\parallel^*$ |
| $t_6$ | | $2 \times (\tau, \tau, \tau)$ <br> $2 \times (1, \tau, \tau)$ | $\dfrac{\tau^3}{12}$ | D | $A_\parallel^*$ |

Construction of icosahedron with pentagonal antiprism and two pentagonal pyramids requires another combination of the tiles $t_5$ and $t_6$ which will be defined as

$$\bar{T}_3 =: t_5 + t_6 + t^5 \tag{9}$$

which is another composite tile with $N_0 = 6$, $N_1 = 10$ and $N_2 = 6$ where the faces consist of Robinson triangles of type $4 \times (1, \tau, \tau)$ and isosceles trapezoids $2 \times (1,1,1,\tau)$. We leave the construction of icosahedron to the reader with the same content of tiles of (7) but with the composite tile $\bar{T}_3$. We will rather follow a composite tiling system with $T_3$ although a tiling system with $\bar{T}_3$ may be possible.

The properties of the fundamental tiles consisting of faces normal to the 5-fold and 3-fold axes can be read from Table 1.

Construction of an icosahedron $i(\tau)$ with edge length $\tau$ is possible with the fundamental tiles. For this we need the inflated tiles of the tiles $t_1, t_2$ and $T_3$ by the inflation factor $\tau$. It is easy to see that $\tau t_1 = t_3 + t_5$ when the face $(1,1,\tau)$ of $t_5$ is matched with a similar face of $t_3$. The inflated $t_2$ can be written as $\tau t_2 = t_4 + t_2 + t_5$ by first matching the equilateral triangular faces of $t_4$ and $t_2$ and then matching the $(1, \tau, \tau)$ faces of $t_2$ and $t_5$. Inflation of the composite tile $T_3$ will be explained in Sec. 3 where $\tau T_3$ is obtained as

$$\tau T_3 = t_1 + 2t_2 + 3t_3 + 4t_4 + 3t_5 + 3t_6. \tag{10}$$

As a result of the inflation of $i(1)$ the icosahedron with edge length $\tau$ will be given as

$$i(\tau) = t_1 + 8t_2 + 10t_3 + 10t_4 + 16t_5 + 3t_6. \tag{11}$$

No further inflation of the icosahedron with the fundamental tiles are possible for the equilateral triangular faces inflated by the factor $\tau$ cannot be dissected into similar triangular faces of edge length 1. But we will see that this is not the case for the dodecahedron.

**b) Projection of the hemicube into a dodecahedron and its dissections**

Let $d(\tau^n)$ denote the dodecahedron of edge length $\tau^n, n = 0, 1, \dots$. A face-first projection of dodecahedron $d(1)$ with its vertices is shown in Fig. 3 where the vertices are taken from the set of vectors (3b). Here $X_i$ and $Y_i$, $i = 1, 2, 3, 4, 5$ are defined by

$$\begin{aligned}
X_1 &= \tfrac{1}{2}(l_1 + l_2 + l_3 + l_4 - l_5 - l_6), \quad X_2 = \tfrac{1}{2}(l_1 + l_2 + l_3 + l_4 + l_5 + l_6), \\
X_3 &= \tfrac{1}{2}(l_1 + l_2 + l_3 - l_4 + l_5 - l_6), \quad X_4 = \tfrac{1}{2}(l_1 - l_2 + l_3 + l_4 + l_5 - l_6) \\
X_5 &= \tfrac{1}{2}(l_1 + l_2 - l_3 + l_4 + l_5 - l_6), \\
Y_1 &= \tfrac{1}{2}(l_1 + l_2 + l_3 - l_4 - l_5 + l_6), \quad Y_2 = \tfrac{1}{2}(l_1 - l_2 + l_3 - l_4 + l_5 + l_6), \\
Y_3 &= \tfrac{1}{2}(l_1 - l_2 - l_3 - l_4 + l_5 - l_6), \quad Y_4 = \tfrac{1}{2}(l_1 - l_2 - l_3 + l_4 - l_5 - l_6), \\
Y_5 &= \tfrac{1}{2}(l_1 + l_2 - l_3 + l_4 - l_5 + l_6).
\end{aligned} \tag{12}$$

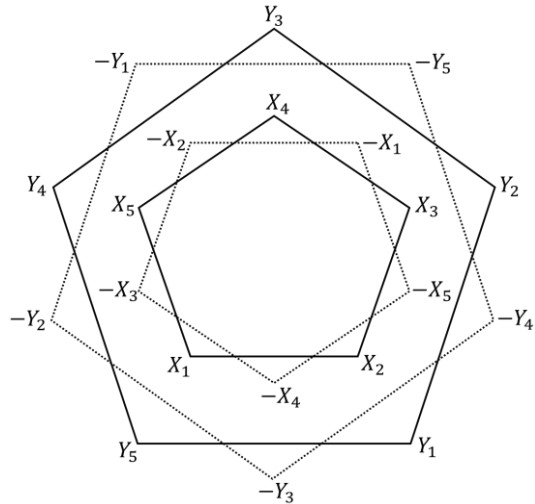

**Figure 3**
Dodecahedron of edge length 1 with vertices $\pm X_i$ and $\pm Y_i$.

In general, the number of a particular facet of a polytope is obtained by Wythoff's technique (Coxeter, 1973, p. 86) by circling a node of the Coxeter-Dynkin diagram. Let us consider the node of the Coxeter- Dynkin diagram of Fig.1 representing the simple root $\alpha_5$ and the corresponding weight vector $\omega_5$: circling a node is equivalent to choosing the weight vector corresponding to this particular simple root. Thus, the orbit of $\omega_5$ under the action of the point group of $D_6$ gives $N_0 = 32$ vertices of the hemicube. To determine the number of vertices ($N_0$), number of edges ($N_1$), number of 2D faces ($N_2$) and in general the number ($N_k$) of the $k$-facet one employs a simple group theoretical technique, namely, the coset decomposition of the point group of $D_6$ under its subgroup(s) fixing the given facet of interest. The number of facets (simplexes in this example) of a hemicube can then be determined as follows

$$\begin{aligned}
N_0 &= \frac{2^5 6!}{6!} = 32, \\
N_1 &= \frac{2^5 6!}{4! 2^2} = 240, \\
N_2 &= \frac{2^5 6!}{3! 3!} = 640, \\
N_3 &= \frac{2^5 6!}{4! 3!} + \frac{2^5 6!}{4! 2} = 180 + 480 = 640, \\
N_4 &= \frac{2^5 6!}{4! 2^4} + \frac{2^5 6!}{5!} = 60 + 192 = 252, \\
N_5 &= \frac{2^5 6!}{5! 2^4} + \frac{2^5 6!}{6!} = 12 + 32 = 44.
\end{aligned} \quad (13)$$

For example, $\omega_5$ is fixed by the generators of the point group $D_6$ forming a subgroup of order 6! so that we obtain $N_0 = 32$ vertices as shown in (13). The rest follows from the subgroups leaving the $k$-facet invariant. Note that they satisfy the Euler characteristic equation $N_0 - N_1 + N_2 - N_3 + N_4 - N_5 = 0$.

This set of formulae can be useful for the projection of the facets of a hemi-cube under the icosahedral group. The factors in the denominators show the orders of the subgroups leaving a given facet invariant which follow from the Coxeter-Dynkin diagram of Fig. 1. Then the number of facets is the sum of the numbers of cosets. First of all, we note that the dodecahedron can be obtained from one of the hemi-cubes which decomposes as the union of a dodecahedron and an icosahedron $32 = 20 + 12$ under the icosahedral group. This implies that the vectors $\pm X_i$ and $\pm Y_i$ represent the vertices of a dodecahedron of edge length 1. The remaining vectors of the hemi-cube define the vertices of an icosahedron of edge length $\tau^{-1}$. Before we proceed further, we remind the reader that each hemi-cube has 640 tetrahedral facets. Each hemi-cube projects into 6 fundamental tiles $t_i, i = 1,2,...,6$ including the tiles $t_1$, $t_2, t_5$ and $t_6$ which constitute the icosahedron $i(1)$ given in (7). We will see that the dodecahedron $d(1)$ can be dissected into the composite tiles defined by

$$\begin{aligned}
T_1 &=: E + C; \ E =: t_4 + t_1 + t_4, \ C =: t_3 + t_6 + t_3, \\
T_2 &=: t_2 + t_4, \\
T_3 &=: t_5 + t_6 + t_5, \\
T_4 &=: t_3 + t_6 + t_5.
\end{aligned} \quad (14)$$

Here the tiles $E$ and $C$ have nothing to the with the notation of Mosseri-Sadoc tiles used for the fundamental tiles. The properties of the composite tiles are displayed in Table 2. Note that $E$ and $C$ are mirror symmetric tiles so as the tile $T_1$. The tile $T_4$ is not mirror symmetric and its mirror reflection can be written as $t_5 + t_6 + t_3$.

**Table 2**
The composite tiles ($N_0$: number of vertices, $N_1$: number of edges, $N_2$: number of faces).

| Name of tile | Figure | $N_0$ | $N_1$ | $N_2$ | Type of faces | Volume |
|---|---|---|---|---|---|---|
| $T_1$ | 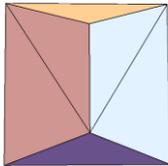 | 8 | 14 | 8 | $4 \times (1,1,\tau)$<br>$4 \times (1,1,1,\tau)$ (trapezoid) | $\dfrac{2\tau^4}{12}$ |
| $T_2$ | 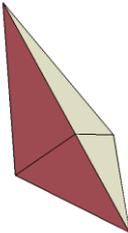 | 4 | 6 | 4 | $2 \times (1,\tau,\tau)$<br>$2 \times (\tau^2,\tau,\tau)$ | $\dfrac{\tau^3}{12}$ |
| $T_3$ | 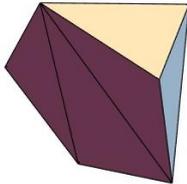 | 6 | 10 | 6 | $5 \times (1,\tau,\tau)$<br>1 pentagon of edge length 1 | $\dfrac{4\tau+3}{12}$ |
| $T_4$ | 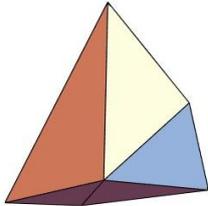 | 6 | 11 | 7 | $3 \times (1,1,\tau)$<br>$3 \times (1,\tau,\tau)$<br>$1(1,1,1,\tau)$ (trapezoid) | $\dfrac{2\tau^3}{12}$ |

The composite tiles are constructed from the fundamental tiles by matching their equilateral triangular faces. This leaves the composite tiles with faces consisting of Robinson triangles only which are all normal to the 5-fold axes. They are vectors which can be represented by the vertices of an icosahedron, say $\pm l_i, (i = 1,2,...,6)$. Since the dual of an icosahedron is a dodecahedron with pentagonal faces which can be constructed with Robinson triangles therefore the normal vectors of the Robinson triangles have 5-fold symmetry. Similarly, equilateral triangles have 3-fold symmetry hence they are orthogonal to the 3-fold symmetric axes of the icosahedral group. The dodecahedra $d(\tau^n)$ of any edge length $\tau^n$ can be constructed with the composite tiles defined by (14) as will be explained in what follows:

The tile $T_1 = E + C$ is made of two composite tiles as it occurs in this combination in any $d(\tau^n)$. The tile $E$ is a nonconvex octahedron obtained by matching two equilateral triangular faces of $t_1$ with equilateral triangular faces of two tiles $t_4$. It has 6 vertices, 12 edges and 8 faces. The set of faces of the composite tile $E$ is $6 \times (1, \tau, \tau)$ and $2 \times (1,1,\tau)$. The tile $C$ is composed by sandwiching a $t_6$ between two $t_3$ tiles at their equilateral triangular faces and has 6 vertices, 12 edges and 8 faces with $6 \times (1,1,\tau)$ and $2 \times (1,\tau,\tau)$. The tile $T_1$ is obtained by inserting $C$ between the legs of $E$ by matching two faces $(1,\tau,\tau)$ with similar faces of $E$. The composite tile $T_1$ now consists of 8 vertices, 14 edges and 8 faces. The tile $T_1$ has $4 \times (1,1,\tau)$ triangles and $4 \times (1,1,1,\tau)$ isosceles trapezoids made of Robinson triangles $(1,1,\tau)$ and $(1,\tau,\tau)$.

The tile $T_2$ which was also used in a 7-tile system of Kramer (Kramer, 1982) is a tetrahedron with faces $2 \times (1,\tau,\tau)$ and $2 \times (\tau^2, \tau, \tau)$ which is obtained by gluing two equilateral faces of tiles $t_2$ and $t_4$. The tile $T_3$ is already described earlier above. The tile $T_4$ is obtained from $T_3$ by replacing one of $t_5$ by $t_3$. Further properties of the composite tiles can be obtained from Table 2. Just to mention another common property is that the dihedral angles between faces of the composite tiles are either $\tan^{-1}(2)$ or $\pi - \tan^{-1}(2)$.

The dodecahedron $d(1)$ can be constructed as the composition of two frustums with $d(1) = d_1(1) + d_2(1)$ where $d_1(1)$ has twice the volume of that of $d_2(1)$. The larger frustum consists of the composite tiles $d_1(1) = 2T_1 + 2T_2 + 3T_4$ and the smaller one $d_2(1) = T_1 + 2T_2 + T_4$ where the composite tiles are assembled face-to-face matching with maximal face coverage. We refer the reader to Appendix A for further details, where dissection of dodecahedron $d(1)$ in terms of the fundamental tiles as well as the composite tiles are studied. These two frustums can be matched at their pentagonal faces of edge length $\tau$ leading to the dodecahedron

$$d(1) = 3T_1 + 4T_2 + 4T_4. \tag{15}$$

The tiles are combined in such a way that they meet at three points defined by $a =: \frac{1}{2}(-\sigma, 1, 0), b =: \frac{1}{2}(-1, 0, -\sigma), c =: \frac{1}{2}(0, \sigma, -1)$ inside the dodecahedron $d(1)$ forming vertices of an equilateral triangle of edge length 1 as can be seen from Appendix A which is just one set of possible assignments of the coordinates for the tiles up to a transformation by the icosahedral group. In fact, by applying the icosahedral transformations one can obtain other assignments of the coordinates. This creates altogether 12 possible intersection points of the tiles in the dodecahedron given by the set of coordinates

$$\tfrac{1}{2}(\pm 1, 0, \pm \sigma), \tfrac{1}{2}(0, \pm \sigma, \pm 1), \tfrac{1}{2}(\pm \sigma, \pm 1, 0), \tag{16}$$

which represent the coordinates of an icosahedron of edge length $\tau^{-1}$. By this we do not mean that an icosahedron of an edge length $\tau^{-1}$ exists in the dodecahedron $d(1)$ rather it admits just three of the coordinates as seen from Appendix A. But this gives an idea that the dodecahedron $d(\tau)$ may embed an icosahedron $i(1)$ of edge length 1. Indeed, this is the case as we will discuss below.

The dodecahedron $d(\tau)$ can be constructed from the icosahedron $i(1)$ with the coordinates of (16) multiplied by $\tau$ by covering the equilateral faces of $t_1$ and $t_2$ with the equilateral faces of $t_4$. By this, one creates a star icosahedron and with this construction the tiles $t_1$ and $t_2$ are converted to the composite tiles $E$ and $T_2$ respectively. Filling the gaps between the legs of $E$ by the tiles $C$ one obtains $7T_1$ composite tiles. The rest follows by face-to-face matching to complete the construction of the dodecahedron given by

$$d(\tau) = 7T_1 + 18T_2 + 14T_3 + 10T_4. \tag{17}$$

Before we proceed further, we should mention that the 12 vertices of the icosahedron inside the dodecahedron $d(\tau)$ exist as intersection points of the composite tiles with no face structures for they are covered by the tiles $t_4$. One can write (17) as the union of inflated tiles as

$$d(\tau) = 3\tau T_1 + 4\tau T_2 + 4\tau T_4. \tag{18}$$

## 3. Composite tiles and the inflation matrix

One can infer an inflation rule with an inflation factor $\tau$ by comparing (17) and (18):

$$\tau E = 2T_2 + T_3 + T_4, \quad \tau C = T_1 + T_3 + T_4, \tag{19}$$

$$\tau T_1 = T_1 + 2T_2 + 2T_3 + 2T_4. \tag{20}$$

Inflation of the other composite tiles can be constructed as follows

$$\begin{aligned}\tau T_2 &= 2T_2 + T_3, \\ \tau T_3 &= T_1 + 2T_2 + T_3 + T_4, \\ \tau T_4 &= T_1 + T_2 + T_3 + T_4,\end{aligned} \tag{21}$$

where $\tau T_i$ ($i = 1, 2, 3, 4$) are illustrated in Fig. 4.

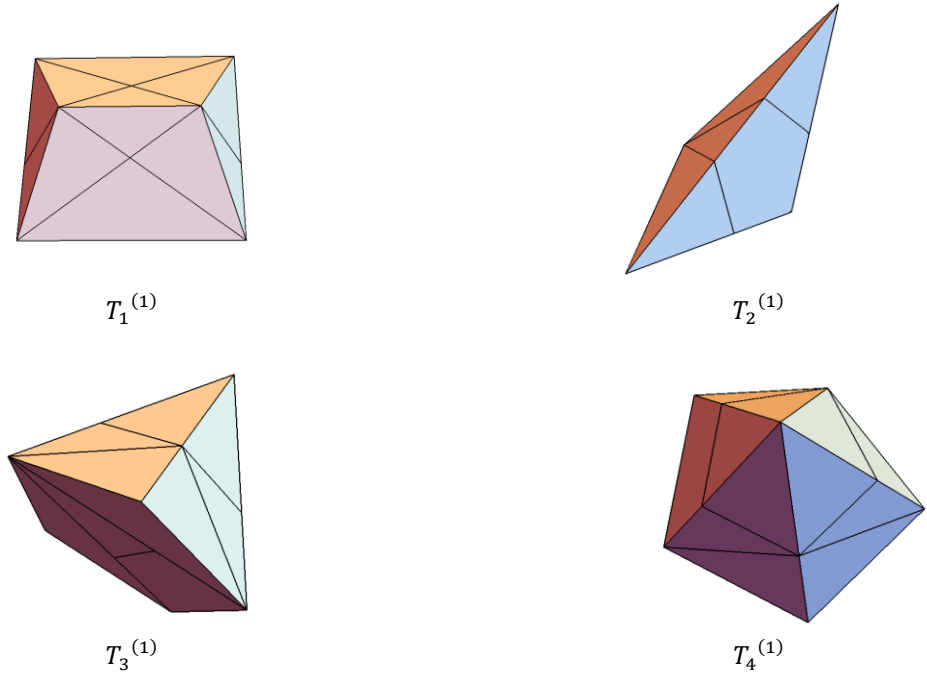

$T_1^{(1)}$  $T_2^{(1)}$

$T_3^{(1)}$  $T_4^{(1)}$

**Figure 4**
The composite tiles $\tau T_i =: T_i^{(1)}$.

The relations (20-21) can be combined in a matrix equation,

$$\tau T_i = \sum_{j=1}^{4} M_{ij} T_j, \ (i = 1, 2, 3, 4),$$

where the matrix $M$ can be written as

$$M = \begin{pmatrix} 1 & 2 & 2 & 2 \\ 0 & 2 & 1 & 0 \\ 1 & 2 & 1 & 1 \\ 1 & 1 & 1 & 1 \end{pmatrix}. \tag{22}$$

The eigenvalues of the inflation matrix are $\tau^3, \tau, \sigma$ and $\sigma^3$; the right eigenvector corresponding to the Perron-Frobenius (PF) eigenvalue $\tau^3$ has the components $(V_{T_1}, V_{T_2}, V_{T_3}, V_{T_4})^T$ with statistical normalization it reads $(\sigma^2, -\frac{\sigma^3}{2}, \frac{4\sigma+3}{2}, -\sigma^3)^T \cong (0.3820, 0.1180, 0.2639, 0.2361)^T$ leading to the relative frequencies of the tiles $T_1, T_2, T_3$ and $T_4$. This implies that the *relative* frequency? of the tile $T_1$ is more than 38%. The statistically normalized left eigenvector or the right eigenvector of $M^T$ of the inflation matrix is $\frac{2}{5\tau+4}(\frac{\tau}{2}, \tau^2, \tau, 1)^T \cong (0.1338, 0.4331, 0.2677, 0.1654)^T$ and it shows the *relative* frequency of the tiles indicating that the tile $T_2$ is nearly 43% more frequent. The PF projection matrix is determined as

$$\lim_{n\to\infty} \tau^{-3n} M^n = P = \frac{1}{30} \begin{pmatrix} 2(\tau+2) & 4(3\tau+1) & 4(\tau+2) & 4\sqrt{5} \\ \sqrt{5} & 2(\tau+2) & 2\sqrt{5} & 2(2+\sigma) \\ 5 & 10\tau & 10 & -10\sigma \\ 2\sqrt{5} & 4(\tau+2) & 4\sqrt{5} & 4(2+\sigma) \end{pmatrix}, \quad P^2 = P. \tag{23}$$

A modification of the tile $T_1$ can be defined as $\acute{T}_1 = T_1 + T_4 = E + C + T_4$ and all equations from (20-23) above can be rephrased with the new set of composite tiles $\acute{T}_1, T_2, T_3$ and $T_4$. For a given tiling system the order of $T_4$ and $C$ in $\acute{T}_1$ does not matter. Equations (15) and (17) can be modified accordingly. The inflation matrix now reads

$$\acute{M} = \begin{pmatrix} 2 & 3 & 3 & 1 \\ 0 & 2 & 1 & 0 \\ 1 & 2 & 1 & 0 \\ 1 & 1 & 1 & 0 \end{pmatrix} \tag{24}$$

The eigenvalues of $\acute{M}$ are again $\tau^3, \tau, \sigma$ and $\sigma^3$; the statistically normalized right eigenvector corresponding to the (PF) eigenvalue $\tau^3$ is $(\frac{1}{2}, \frac{\sigma^2}{4}, \frac{-3\sigma-1}{4}, \frac{\sigma^2}{2})^T$ which implies that the *relative* frequency of the tile $\acute{T}_1$ is exactly 50%. Statistically normalized right eigenvector of transpose of the matrix $\acute{M}$ now reads $(\frac{1}{4\tau}, \frac{1}{2}, \frac{1}{2\tau}, \frac{1}{4\tau^4})^T$. Here again the *relative* frequency of the tile $T_2$ is 50%. The equations (15) and (17) now read respectively

$$d(1) = 3\acute{T}_1 + 4T_2 + T_4 = 3(\acute{T}_1 + T_2) + (T_2 + T_4), \tag{25}$$

and

$$d(\tau) = 7\acute{T}_1 + 18T_2 + 14T_3 + 3T_4. \tag{26}$$

As we will show in Appendix A that there is a 3-fold symmetric construction of $d(1)$ so that three sets of $(\acute{T}_1 + T_2)$ are permuted into each other while $(T_2 + T_4)$ remains intact. This is a result of the triangular symmetry of the vertices $(a, b, c)$ as pointed out above.

Let us continue with the original set of composite Mosseri-Sadoc tiles although similar studies can be carried out with the new set of tiles.

After this general procedure we will illustrate some of the inflated patches. For this, we first define the inflated tiles by a new notation. Let us denote by $\tau^n T_i =: T_i^{(n)}, n = 0, 1, 2, ...$ then we can write

$$T_i^{(n)} = \sum_{j=1}^{4} (M^n)_{ij} T_j, \quad (i = 1, 2, 3, 4). \tag{27}$$

The dodecahedron $d(\tau^n)$, $(n = 0, 1, 2, ...)$ can be constructed in terms of the composite tiles and they can also be dissected in terms of the dodecahedra $d(1)$ and $d(\tau)$ along with the composite tiles. For example, we obtain the dodecahedron $d(1)$ in the inflated tiles given by

$$\begin{aligned} T_1^{(2)} &= d(1) + 2T_2^{(1)} + T_3^{(1)} + T_4^{(1)} + T_2 + 4T_3, \\ T_2^{(3)} &= d(1) + 2T_2^{(2)} + 5T_2 + 6T_3, \\ T_3^{(2)} &= d(1) + 5T_2 + 6T_3, \\ T_4^{(2)} &= d(1) + 3T_2 + 5T_3, \end{aligned} \quad (28)$$

where they are depicted in Fig. 5.

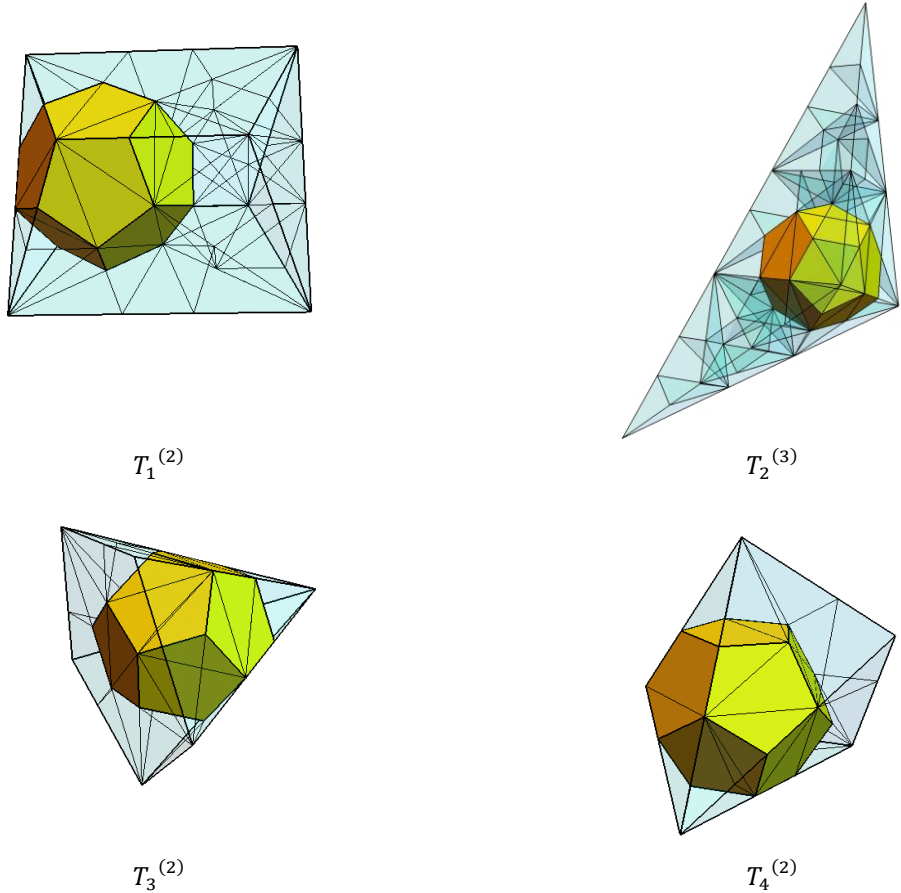

$T_1^{(2)}$    $T_2^{(3)}$

$T_3^{(2)}$    $T_4^{(2)}$

**Figure 5**
The inflated composite tiles of (28) with dodecahedral structures.

Further inflation of tiles in (28) by $\tau$ will produce $d(\tau)$ along with $d(1)$, for example, in the simpler case we obtain

$$T_2^{(4)} = 2d(1) + d(\tau) + 4T_2^{(2)} + 5T_2^{(1)} + 6T_3^{(1)} + 10T_2 + 12T_3. \quad (29)$$

Another interesting case happens in the inflation represented by

$$T_1^{(4)} = 13d(1) + 2d(\tau) + 9T_2^{(2)} + 14T_2^{(1)} + 14T_3^{(1)} + 3T_4^{(1)} + 45T_2 + 68T_3, \quad (30)$$

where dodecahedra $13d(1)$ and $2d(\tau)$ occur simultaneously. Similar formulae can be obtained for $T_3^{(4)}$ and $T_4^{(4)}$ as

$$T_3^{(4)} = 13d(1) + 9T_2^{(2)} + 6T_2^{(1)} + 3\,T_3^{(1)} + 3T_4^{(1)} + 45\,T_2 + 68T_3, \tag{31}$$

$$T_4^{(4)} = 12d(1) + 7T_2^{(2)} + 6T_2^{(1)} + 3\,T_3^{(1)} + 3T_4^{(1)} + 40\,T_2 + 62T_3. \tag{32}$$

They are illustrated in Fig. 6 by highlighting the dodecahedral structures and leaving the others transparent.

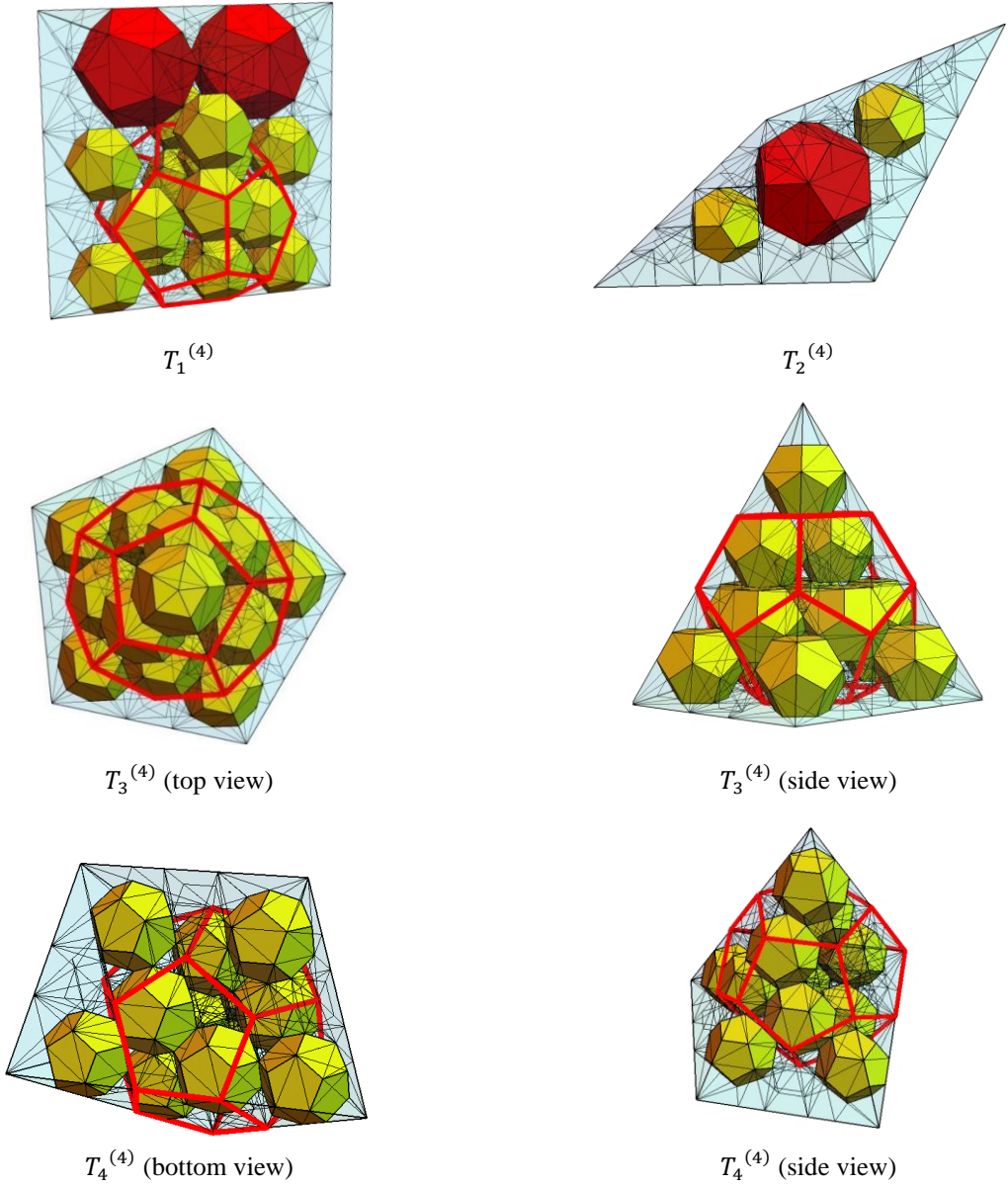

$T_1^{(4)}$ \qquad $T_2^{(4)}$

$T_3^{(4)}$ (top view) \qquad $T_3^{(4)}$ (side view)

$T_4^{(4)}$ (bottom view) \qquad $T_4^{(4)}$ (side view)

**Figure 6**
Dodecahedral structures in $T_i^{(4)}$, $i = 1, 2, 3, 4$ where the tiles $d(1)$ and $d(\tau)$ are demonstrated in different colours and the dodecahedral frames is denoted by $d(\tau^2)$.

The dodecahedron $d(1)$ can be inflated to an arbitrary order of the inflation factor $\tau$. They all reduce to a number of $d(1)$ and/or $d(\tau)$ along with the other composite tiles. As we have discussed in Appendix A and Appendix B it is possible to dissect the dodecahedron to obtain 2-fold, 3-fold and 5-fold symmetric distributions of dodecahedra. To give an example we illustrate $d(\tau^2)$ in Fig. 7 which consists of $7d(1)$ and accompanying composite tiles left transparent.

Equation (25) proves that $7d(1)$ decomposes as $7d(1) = 3d(1) + 3d(1) + d(1)$ implying that one of the dodecahedron is invariant under the dihedral group of order 6 while the rest decompose into two sets of three-fold symmetric combinations which is discussed in Appendix A. This construction of dodecahedron reveals all symmetries of dodecahedron. In Appendix B we give an alternative dissection of dodecahedron which displays 5-fold symmetry.

It is clear that the patches include $d(1)$ and $d(\tau)$ in abundance. We illustrate some patches from $d(\tau^3)$, $d(\tau^4)$ and $d(\tau^5)$ displaying 5-fold, 3-fold and 2-fold symmetries around certain axes as depicted in Fig. 8 with transparent composite tiles. The dodecahedra of interest are given by

$$d(\tau^3) = 10d(1) + 7d(\tau) + 46\, T_1 + 222 T_2 + 146 T_3 + 46 T_4, \qquad (33)$$

$$d(\tau^4) = 95d(1) + 10d(\tau) + 170\, T_1 + 1110 T_2 + 898 T_3 + 170 T_4, \qquad (34)$$

$$d(\tau^5) = 240d(1) + 95d(\tau) + 828\, T_1 + 4446 T_2 + 3078 T_3 + 828 T_4. \qquad (35)$$

The set of dodecahedra in (33-35) has 5-fold, 3-fold and 2-fold symmetry axes where blue and gold colors represent the dodecahedra $d(\tau)$ and $d(1)$ respectively in $d(\tau^3)$. The set of dodecahedra of $d(\tau^4)$ is displayed viewing it from three symmetry axes clustered in $d(1)$(gold), $d(\tau)$(blue) and $d(\tau^2)$(red). Similarly, the dodecahedron $d(\tau^5)$ has been displayed from 3-fold and 5-fold symmetry axes first with blue and gold dodecahedra and then their clusters are combined in different colors of $d(1), d(\tau), d(\tau^2)$(light blue) and $d(\tau^3)$(red).

The numbers of dodecahedra $d(1)$ and $d(\tau)$ in a dodecahedron $d(\tau^n)$ grow fast as the order of inflation by $\tau$ inreases.

To give an example consider the dodecahedron $d(1)$ inflated by $\tau^{10}$,

$$d(\tau^{10}) = 432139 d(1) + 92850 d(\tau) + 1064050 T_1 + 6341550 T_2 + 4720730 T_3 + 1064050 T_4. \quad (36)$$

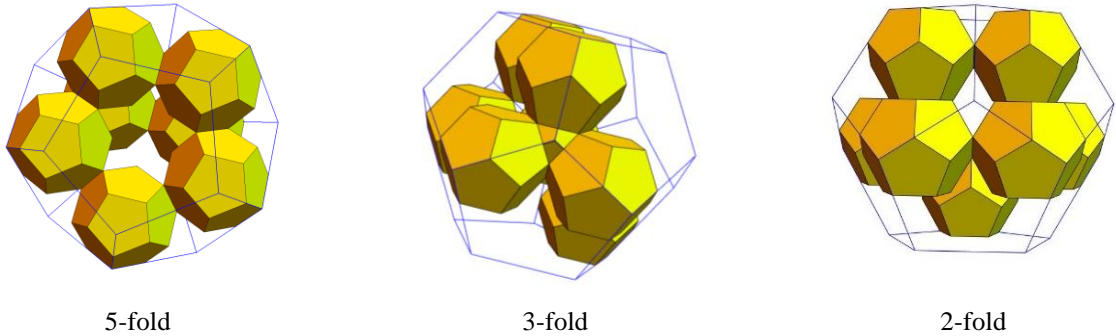

| 5-fold | 3-fold | 2-fold |

**Figure 7**
Symmetries of $d(\tau^2)$.

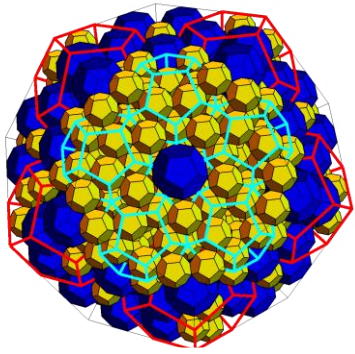
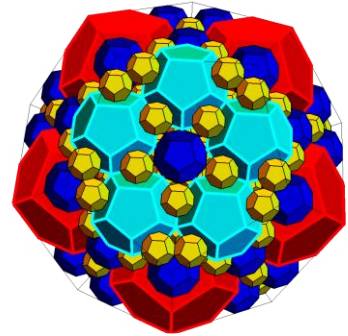

$d(\tau^5)$: 5-fold symmetry

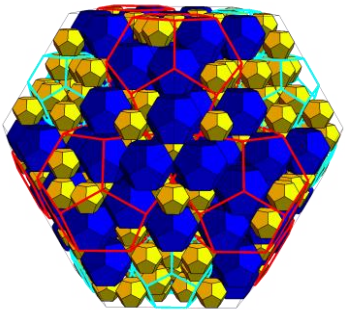
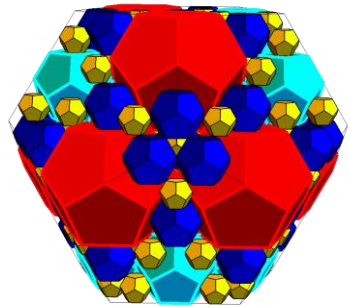

$d(\tau^5)$: 3-fold symmetry

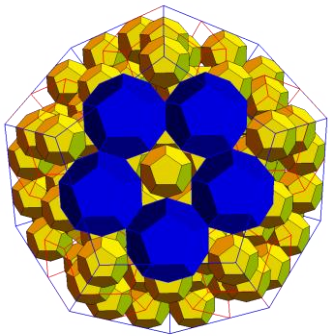
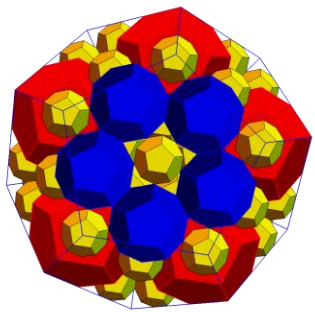
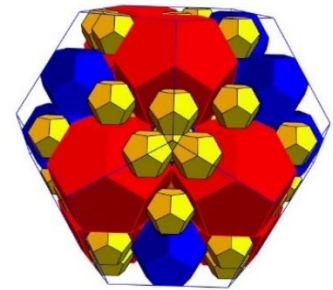

$d(\tau^4)$: 5-fold, 5-fold, 3-fold symmetries

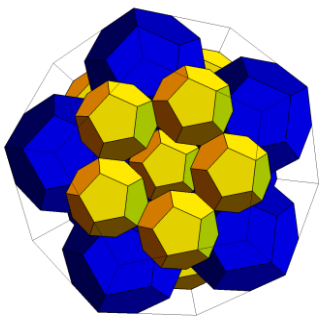
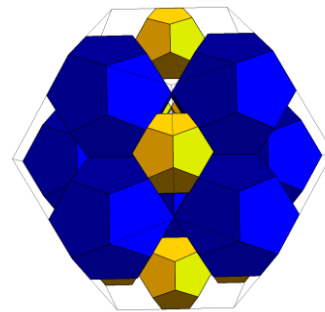

$d(\tau^3)$: 5-fold symmetry                    $d(\tau^3)$: 2-fold symmetry

**Figure 8.** Patches of dodecahedra $d(\tau^5)$, $d(\tau^4)$ and $d(\tau^3)$ with 5-fold, 3-fold and 2-fold symmetries (for further details see the text).

## 4. Concluding remarks

The present tiling scheme with its dodecahedral structures does not only show that the 3D-Euclidean space can be tiled with icosahedral symmetry but it may prove that it is an alternative model to the icosahedral quasicrystals since some experiments display dodecahedral structures such as the holmium-magnesium-zinc quasicrystal (Tsai et. al., 1994). The three-fold embedding of the composite tiles in the dodecahedron may have far reaching consequences. We have exhibited 5-fold, 3-fold and 2-fold symmetries of the icosahedral group as clusters of dodecahedra. We have used a face-to-face tiling scheme with maximal face coverage. Faces of the inflated tiles exhibit planar tilings with Robinson triangles. The Fundamental tiles allow more composite tiles such as $\bar{T}_3 =: t_5 + t_6 + t^5$ which is used in the construction of icosahedron which may lead to different tiling models. We should also emphasize that icosahedron and/or icosidodecahedron are dissected in terms a mixture of fundamental and composite tiles and cannot be inflated beyond the edge lengths 1 and $\tau$.

## Appendix A. Dissection of dodecahedron $d(1) = 3(T_1 + T_4 + T_2) + (T_2 + T_4)$

$$a =: \frac{1}{2}(-\sigma, 1, 0), b =: \frac{1}{2}(-1, 0, -\sigma), c =: \frac{1}{2}(0, \sigma, -1)$$

### *Embedding of $T_1 + T_4 + T_2$ in dodecahedron with three-fold symmetry*

The dodecahedron can be dissected as the 3-fold symmetric combination of the set of tiles $T_1 + T_4 + T_2$. Below we illustrate coordinates of this symmetric combinations.

$$E + C + T_4$$

$$\begin{array}{lll}
t_4: (X_5, X_1, a, -Y_2); & t_3: (-Y_2, -Y_1, Y_4, X_5); & t_3: (Y_3, X_4, -Y_1, -Y_5) \\
t_1: (X_5, X_1, a, X_4); & t_6: (-Y_2, -Y_1, X_5, a); & t_6: (a, X_4, -Y_1, -Y_5) \\
t_4: (X_5, X_4, a, -Y_1); & t_3: (-Y_2, -Y_1, a, c); & t_5: (c, a, -Y_1, -Y_5)
\end{array}$$

$$T_2$$

$$t_2: (X_1, X_3, X_4, a); t_4: (X_3, X_4, a, -Y_5). \tag{A1}$$

The other sets of tiles of $E + C + T_4$ and $T_2$ can be obtained by the cyclic permutation of the coordinates as described below. The cyclic transformation between the coordinates $(a, b, c)$ leaves the vertices $\pm Y_5$ invariant and induces the transformations between the sets of pentagonal vertices

$$\pm (X_1, X_2, X_3, X_4, X_5) \to \pm(-Y_3, -X_4, -X_5, -Y_4, Y_1) \to \pm(-Y_2, Y_4, -Y_1, -X_2, -X_3), \tag{A2}$$

$$\pm (Y_2, -Y_4, -X_5, -X_1, -Y_5) \to \pm(-X_1, -X_2, -Y_1, Y_3, -Y_5) \to \pm(Y_3, X_4, X_3, Y_2, -Y_5). \tag{A3}$$

### *Embedding of $T_2 + T_4$ in dodecahedron*

The tile $(T_2 + T_4)$ is invariant under the 3-fold symmetry. It can be represented by three sets of coordinates however all coordinates describe the same union $(T_2 + T_4)$. For this, it suffices to give one set of vertices only

$$T_2: \{t_2: (a, b, c, -Y_3), t_4: (a, b, c, -Y_5)\}, \tag{A4}$$

$$T_4: \{t_3: (-Y_3, -Y_2, X_1, Y_5); t_6: (X_1, a, -Y_2, -Y_3); t_5: (-Y_2, -Y_3, c, a)\}. \tag{A5}$$

Volume of $d(1)$ is given by

$$\text{Vol}(d(1)) = \frac{1}{12}(6\tau^4 + 4\tau^3 + 8\tau^3) = \frac{1}{2}(7\tau + 4). \tag{A6}$$

## Appendix B. Dissection of $d(1) = d_1(1) + d_2(1)$ leading to five-fold symmetry

Here we give the coordinates of the tiles in $d(1)$ leading to the 5-fold symmetric dodecahedral distribution. The tiles $C$ and $T_4$ generate dodecahedra at the second order of inflation. The large frustum $d_1(1)$ consists of $2C$ and $3T_4$ so that we observe 5 dodecahedra in the second order of inflation leading to 5-fold dodecahedral symmetry.

**Coordinates of $d_1(1) = 2T_1 + 2T_2 + 3T_4$**

$$\{t_1: [X_1, X_4, X_5, a]; t_4: [X_1, X_5, -Y_2, a]; t_4: [X_4, X_5, -Y_1, a]\} = E,$$

$$\{t_6: [X_5, -Y_1, -Y_2\, a]; t_3: [X_5, -Y_1, -Y_2, Y_4]; t_3: [-Y_1, -Y_2, a, c]\} = C,$$

$$T_1 = E + C \tag{B1}$$

$$\{t_1: [X_2, X_3, X_4, a]; t_4: [X_3, X_4, -Y_5, a]; t_4: [X_2, X_3, -Y_4, a]\} = E,$$

$$\{t_6: [X_3, -Y_4, -Y_5, a]; t_3: [X_3, -Y_4, Y_2, -Y_5]; t_3: [-Y_4, -Y_5, c, a]\} = C,$$

$$T_1 = E + C \tag{B2}$$

$$\{t_2: [X_1, X_2, X_4, a]; t_4: [X_1, X_2, -Y_3, a\} = T_2 \tag{B3}$$

$$\{t_4: [-Y_5, b, a, c]; t_2: [-Y_3, b, a, c\} = T_2 \tag{B4}$$

$$\{t_3: [X_1, -Y_2, Y_5, -Y_3]; t_6: [X_1, -Y_2, -Y_3, a]; t_5: [-Y_2, -Y_3, a, c]\} = T_4 \tag{B5}$$

$$\{t_3: [-Y_3, Y_1, -Y_4, X_2]; t_6: [X_2, -Y_3, -Y_4, a]; t_5: [-Y_3, -Y_4, a, b]\} = T_4 \tag{B6}$$

$$\{t_3: [X_4, Y_3, -Y_1, -Y_5]; t_6: [X_4, -Y_1, -Y_5, a]; t_5: [-Y_1, -Y_5, a, c]\} = T_4 \tag{B7}$$

**Coordinates of $d_2(1) = T_1 + 2T_2 + T_4$**

$$\{t_1: [-X_4, -Y_4, -X_5, b]; t_4: [-X_4, -X_5, -X_2, b]; t_4: [-Y_4, -X_5, -Y_5, b\,]\} = E,$$

$$\{t_6: [-X_5, -X_2, -Y_5, b]; t_3: [-X_1, -X_2, -X_5, -Y_5]; t_3: [-X_2, -Y_5, b, c]\} = C,$$

$$T_1 = E + C \tag{B8}$$

$$\{t_2: [-Y_2, -Y_1, -X_2, c]; \ t_4: [-X_2, -Y_1, -Y_5, c]\} = T_2 \quad (B9)$$

$$\{t_2: [-Y_3, -Y_4, -X_4, b]; \ t_4: [-X_4, -Y_2, -Y_3, b]\} = T_2 \quad (B10)$$

$$\{t_3: [-X_4, -X_3, -X_2, -Y_2]; \ t_6: [-X_4, -X_2, -Y_2, b]; \ t_5: [-X_2, -Y_2, b, c]\} = T_4 \quad (B11)$$

The 12 pentagonal faces of dodecahedron are represented by the sets of 5 vertices,

$$\pm(X_1, X_2, X_3, X_4, X_5); \ \{(X_1, X_2, Y_1, -Y_3, Y_5) \text{ and cyclic permutations of } X_i \text{ and } Y_i\}. \quad (B12)$$

**Appendix C. Construction of icosidodecahedron in terms of the fundamental tiles**

The root system of $D_6$ is obtained from the point group application on the weight vector $\omega_2 = l_1 + l_2$ leading to the vertices of the root polytope

$$\pm l_i \pm l_j, \quad i \neq j = 1, 2, \ldots, 6. \quad (C1)$$

Its projection into 3D-space is the union of two icosidodecahedra of edge lengths 1 and $\tau^{-1}$ where $l_1 + l_2$ generates one orbit of 30 vertices of icosidodecahedron of edge length 1 and the vertex $l_1 - l_2$ generates another orbit of edge length $\tau^{-1}$. An icosidodecahedron consists of 30 vertices, 32 faces (12 pentagonal and 20 triangular faces) and 60 edges. A symmetric half of the vertices of the icosidodecahedron of edge length 1 is depicted in Fig. 9 where the other half represents the negatives of those in Fig. 9. The 15 vertices of the icosidodecahedron can be obtained from the 5-fold symmetry $(1)(235\bar{6}4)$ and the 3-fold symmetry $(123)(465)$.

The 30 vertices of icosidodecahedron can be written as

$$\tau\{(\pm 1, 0, 0), (0, \pm 1, 0), (0, 0, \pm 1), \tfrac{1}{2}(\pm \tau \pm \sigma \pm 1), \tfrac{1}{2}(\pm \sigma \pm 1 \pm \tau), \tfrac{1}{2}(\pm 1 \pm \tau \pm \sigma)\}. \quad (C2)$$

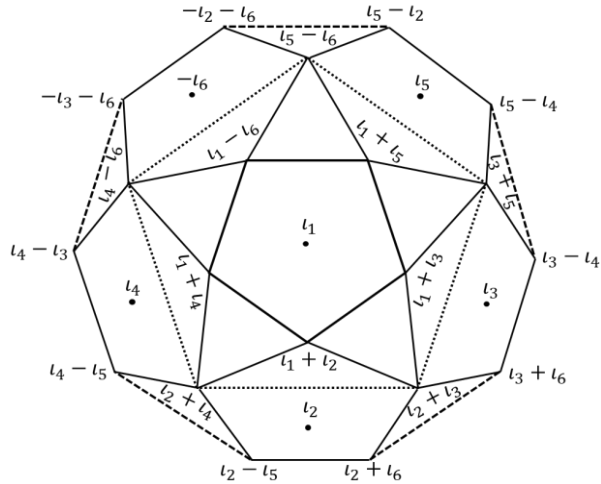

**Figure 9**
A sketch of half icosidodecahedron. (The vectors $\pm l_i$ don't represent vertices but illustrate the vectors normal to the pentagonal faces).

Half of the icosidodecahedron consists of 6 pentagonal and 10 triangular faces; it is then straightforward to see that 6 tiles of $T_3 =: t_5 + t_6 + t_5$ and 10 fundamental tiles of type $t_4$ fill half the icosidodecahedron where the pentagonal bases of $T_3$ and the triangular faces of $t_4$ coincide with the corresponding faces of the icosidodecahedron. The apexes lie at the origin of the coordinate system. Two halves of icosidodecahedron are brought face-to-face to complete the icosidodecahedral structure denoted by $id(1)$ which now consists of 12 $T_3$ and 20 tiles of $t_4$ or in terms of fundamental tiles

$$id(1) = 20t_4 + 12T_3 = 20t_4 + 24t_5 + 12t_6. \tag{C3}$$

Volume of an icosidodecahedron $id(1)$ is the sum of the volumes of its constituents and is given as

$$\text{Vol}(id(1)) = \tfrac{1}{12}(20\tau^2 + 24\tau^2 + 12\tau^3) = \tfrac{1}{3}(17\tau + 14). \tag{C4}$$

The icosidodecahedron $id(\tau)$ of edge length $\tau$ can be obtained from the inflation of $t_4$ and $T_3$ by the inflation factor $\tau$:

$$id(\tau) = 20\tau t_4 + 12\tau T_3. \tag{C5}$$

One can easily see that an inflation of $t_4$ by $\tau$ is given by

$$\tau t_4 = \tau t_2 + t_6 = t_4 + t_2 + t_5 + t_6, \tag{C6}$$

where $t_4 + t_2$ are matched on their equilateral faces of edge length 1, $t_2 + t_5$ on their faces of Robinson triangles $(1, \tau, \tau)$ and $t_5 + t_6$ on their equilateral faces of edge length $\tau$ respectively. The inflation of $T_3$ is already given by (21) so that dissection of icosidodecahedron of edge length $\tau$ in terms of the fundamental tiles reads

$$id(\tau) = 12t_1 + 44t_2 + 36t_3 + 68t_4 + 56t_5 + 56t_6. \tag{C7}$$